\documentclass[12pt]{amsart}
\usepackage{latexsym, amssymb, amscd, rotating}
\renewcommand\a{\alpha}
\renewcommand\b{\beta}
\newcommand\g{\gamma}
\renewcommand\d{\delta}
\newcommand\la{\lambda}
\newcommand\z{\zeta}

\newcommand\m{\mu}

\newcommand\vf{\varphi}

\renewcommand\r{\rho}

\newcommand\Om{\Omega}

\newcommand\vL{\varLambda}

\newcommand\ve{\varepsilon}

\newcommand{\QQ}{\mathbb Q}
\newcommand{\FF}{\mathbb F}

\newcommand{\ZZ}{\mathbb Z}

\newcommand{\CC}{\mathbb C}

\newcommand\Bi{\mathbf i}

\newcommand\Bm{\mathbf m}
\newcommand\Bs{\mathbf s}

\newcommand\MC{\mathcal{C}}

\newcommand\CI{\mathcal{I}}

\newcommand\CM{\mathcal{M}}

\newcommand\CP{\mathcal{P}}

\newcommand\CF{\mathcal{F}}

\newcommand\CX{ \mathcal{X}}

\newcommand\FS{\mathfrak S}

\newcommand\Bdel{\boldsymbol\delta}
\newcommand\BLa{\boldsymbol\varLambda}

\newcommand\Bmu{\boldsymbol\mu}

\newcommand\Ba{\boldsymbol\alpha}
\newcommand\Bb{\boldsymbol\beta}
\newcommand\Bg{\boldsymbol\gamma}

\newcommand\iv{^{-1}}

\newcommand\wt{\widetilde}

\newcommand\ol{\overline}

\newcommand\id{\operatorname{id}}
\newcommand\lp{\operatorname{\langle\!}}
\newcommand\rp{\operatorname{\!\rangle}}

\newcommand{\isom}{\,\raise2pt\hbox{$\underrightarrow{\sim}$}\,}
\numberwithin{equation}{section}

\newtheorem{thm}{Theorem}[section]
\newtheorem{lem}[thm]{Lemma}

\newtheorem{prop}[thm]{Proposition}

\def \para#1{\par\medskip\textbf{#1}
              \addtocounter{thm}{1}}
 
\begin{document}
\setlength{\baselineskip}{4.9mm}
\setlength{\abovedisplayskip}{4.5mm}
\setlength{\belowdisplayskip}{4.5mm}
\renewcommand{\theenumi}{\roman{enumi}}
\renewcommand{\labelenumi}{(\theenumi)}
\renewcommand{\thefootnote}{\fnsymbol{footnote}}
\parindent=20pt
\medskip
\title{Macdonald functions associated to 
      \\  complex reflection groups}
\author{Toshiaki Shoji}
\maketitle
\pagestyle{myheadings}
\markboth{T. SHOJI}{MACDONALD FUNCTIONS}
\begin{center}
Department of Mathematics  \\
Science University of Tokyo \\ 
Noda, Chiba 278-8510, Japan
\end{center}
\par\medskip
\begin{center}
{\it To Robert Steinberg}
\end{center}
\begin{abstract}
Let $W$ be the complex reflection group 
$\FS_n\ltimes (\ZZ/e\ZZ)^n$.
In the author's previous paper [S1], Hall-Littlewood functions 
associated to $W$ were introduced.  In the special case
where $W$ is a Weyl group of type $B_n$, they are closely
related to Green polynomials of finite classical groups.
In this paper, we introduce a two variables version of 
the above Hall-Littlewood functions, as a generalization 
of Macdonald functions associated to symmetric groups.  
A generalization of Macdonald operators is also constructed,
and we characterize such functions 
by making use of Macdonald operators, assuming a certain 
conjecture.  
\end{abstract}
\bigskip
\medskip
\addtocounter{section}{-1}
\section{Introduction}
Macdonald functions $P_{\la}(x;q,t)$, which  were introduced by
I.G. Macdonald [M2] in 1987, are two variables
versions of Hall-Littlewood functions $P_{\la}(x;t)$.
Those Hall-Littlewood functions and Macdonald functions are
parametrized by partitions $\la$ of $n$.  Since partitions of 
$n$ parameterize irreducible characters of the symmetric group
$\FS_n$, we may say that these functions are associated to
symmetric groups.  On the other hand, Hall-Littlewood functions
are closely related to Green polynomials of a finite general linear
group $GL_n(\FF_q)$.  In this direction, partitions $\la$ of $n$ 
occur as unipotent classes of $GL_n(\FF_q)$.  
Unipotent classes of other finite classical groups such as 
$Sp_{2n}(\FF_q), SO_{2n+1}(\FF_q)$ have more complicated patterns.
Lusztig introduced in [L2] (unipotent) symbols, as a generalization of the
notion of partitions, to describe such unipotent classes in connection
with Springer representations of Weyl groups.  He also introduced in
[L1] a notion of symbols to parameterize unipotent characters of 
finite classical groups.  
\par  
In [S1], the author constructed Hall-Littlewood functions 
associated to complex reflection groups 
$W \simeq \FS_n\ltimes (\ZZ/e\ZZ)^n$.  In [S2], [S3],
some related topics are discussed.  
Our Hall-Littlewood functions are parametrized
by $e$-tuples of partitions (which parameterize irreducible characters of
$W$), or rather by various types of $e$-symbols.  In the case where
$e=2$, $W$ is the Weyl group of type $B_n$.  In this case, 
Hall-Littlewood functions attached to unipotent symbols are closely
related to Green functions of $Sp_{2n}(\FF_q)$ or $SO_{2n+1}(\FF_q)$.
It is also expected that our Hall-Littlewood functions attached to
symbols have some connection with unipotent characters.
\par
This paper is an attempt to generalize Hall-Littlewood functions 
$P^{\pm}_{\BLa}(x;t)$ (where $\BLa$ is an $e$-symbol) to the two variables
version $P^{\pm}_{\BLa}(x;q,t)$, just as in the case of
$P_{\la}(x;q,t)$.  We call such functions $P^{\pm}_{\BLa}(x;q,t)$
Macdonald functions associated to  $W$. 
In the case of original Macdonald functions, they are characterized 
as simultaneous eigenfunctions of various Macdonald operators.
We also construct Macdonald operators having  (conjecturally) good
properties. However, we note that our construction of Macdonald 
operators works only for a special type of symbols (in the case
where $e = 2$, this is exactly the symbols used to parameterize
unipotent characters), though Macdonald functions can be
constructed for any type of symbols.  In the case of original 
Macdonald operators,
the representation matrix with respect to the basis of Schur functions
is a triangular matrix, with distinct eigenvalues.  In our case,
the matrix of Macdonald operators turns out to be a block triangular
matrix, where the blocks correspond to families of symbols.
One can conjecture that the matrices appearing in the diagonal blocks
have no common eigenvalues.  Assuming this conjecture, we show 
that the Macdonald operator characterizes Macdonald functions, 
not as eigenfunction, but as a unique solution of linear systems 
attached to the above diagonal blocks.    
\par
The properties of Macdonald functions discussed in this paper 
are just a part of those established for the original 
Macdonald functions.  We hope to discuss more about them 
in a subsequent paper.     
    
\par\bigskip
\section{Symmetric functions with two parameters}
\para{1.1.}
An $e$-tuple of partitions $\Ba = (\a^{(0)}, \dots, \a^{(e-1)})$ is 
called an $e$-partition.  We define the size $|\Ba|$ of $\Ba$ by  
$|\Ba| = \sum_{k =0}^{e-1} |\a^{(k)}|$, where 
$|\a^{(k)}|$ is the size of the partition $\a^{(k)}$.
For a partition $\a: \a_1 \ge \a_2 \ge \cdots \ge \a_k\ge 0$, 
let $l(\a)$ be the number of non-zero parts $\a_k$.
We denote by $\CP_n$ the set of partitions of $n$, and 
$\CP_{n,e}$ the set of $e$-partitions of size $n$.  
Let $W$ be the complex reflection group $\FS_n\ltimes (\ZZ/e\ZZ)^n$.
Then the set of conjugacy classes in $W$ is in one to one correspondence
with the set $\CP_{n,e}$.
\par
Let us fix a sequence of positive integers 
$\Bm = (m_0, \dots, m_{e-1})$, and consider indeterminates 
$x_j^{(k)} \ (0 \le k \le e-1, 1 \le j \le m_k)$.
We denote by $x = x_{\Bm}$ the whole variables $(x_j^{(k)})$, and 
also denote by $x^{(k)}$ the variables 
$x_1^{(k)}, \dots, x_{m_k}^{(k)}$.
Let $\z$ be a primitive $e$-th root of unity in
$\CC$.
For each integer $r \ge 1$ and $i$ such that $0 \le i \le e-1$, put
\begin{equation*}
p_r^{(i)}(x) = \sum_{j=0}^{e-1}\z^{ij}p_r(x^{(j)}),
\end{equation*}
where $p_r(x^{(j)})$ denotes the $r$-th power sum symmetric function
with respect to the variables $x^{(j)}$.
We put $p_r^{(i)}(x) = 1$ for $r = 0$.
For an $e$-partition $\Ba = (\a^{(0)}, \dots, \a^{(e-1)})$ with 
$\a^{(k)}: \a^{(k)}_1 \ge \dots \ge \a^{(k)}_{m_k}$, we
define a function $p_{\Ba}(x)$ by
\begin{equation*}
p_{\Ba}(x) = \prod_{k=0}^{e-1}\prod_{j=1}^{m_k}p^{(k)}_{\a^{(k)}_j}(x).
\end{equation*}
\par
Next, we define the Schur function $s_{\Ba}(x)$ and monomial symmetric
functions $m_{\Ba}(x)$ associated to $\Ba$ by
\begin{equation*}
s_{\Ba}(x) = \prod_{k=0}^{e-1}s_{\a^{(k)}}(x^{(k)}), \qquad
m_{\Ba}(x) = \prod_{k=0}^{e-1}m_{\a^{(k)}}(x^{(k)}), 
\end{equation*}
where $s_{\a^{(k)}}(x^{(k)})$ (resp. $m_{\a^{(k)}}(x^{(k)}$ ) 
denotes the usual Schur function (resp. monomial symmetric function)
associated
to the partition $\a^{(k)}$ with respect to the variables 
$x^{(k)}$.
These are the symmetric functions associated to complex 
reflection groups $W \simeq \FS_n\ltimes (\ZZ/e\ZZ)^n$, as given 
in [M1, Appendix B].
\para{1.2.}
Put $\FS_{\Bm} = \FS_{m_0}\times \cdots \times \FS_{m_{e-1}}$.
We denote by  
$\Xi_{\Bm} = \bigotimes_{k = 0}^{e-1}
       \ZZ[x_1^{(k)}, \dots, x_{m_k}^{(k)}]^{\FS_{m_k}}$ 
the ring of symmetric polynomials (with respect to $\FS_{\Bm}$) 
with variables $x = (x_j^{(k)})$.
$\Xi_{\Bm}$ has a structure of a graded ring 
$\Xi_{\Bm} = \bigoplus_{i \ge 0}\Xi_{\Bm}^i$, where $\Xi_{\Bm}^i$ consists of
homogeneous symmetric polynomials of degree $i$, 
together with the zero polynomial.  We consider the inverse limit 
$$
\Xi^i = \lim_{\substack{\leftarrow \\ {\Bm}}}\ \Xi_{\Bm}^i
$$
with respect to homomorphisms 
$\r_{\Bm', \Bm} : \Xi_{\Bm'}^i \to \Xi_{\Bm}^i$, 
where $\Bm' = (m_0', \dots,m_{e-1}')$ with $m_k' = m_k + l$ for some 
integer $l \ge 0$, and $\r_{\Bm',\Bm}$ is induced from  
the homomorphism 
$\bigotimes_k\ZZ[x_1^{(k)}, \dots, x_{m'_k}^{(k)}] \to 
\bigotimes_k\ZZ[x_1^{(k)}, \dots, x_{m_k}^{(k)}] $  given by
sending $x_i^{(k)}$ to 0 for $i > m_k$, and leaving the other
$x_j^{(k)}$ invariant.   $\Xi = \bigoplus_{i \ge 0}\Xi^i$ is called 
the space of symmetric functions.  Schur functions $s_{\Ba}(x)$ with
infinitely many variables $x^{(k)}_1, x^{(k)}_2 \dots$ are regarded as 
elements in $\Xi^n$ with $n = |\Ba|$, and the set 
$\{ s_{\Ba}(x) \mid \Ba \in \CP_{n,e}\}$ forms a $\ZZ$-basis of $\Xi^n$.  
Similarly, $\{ m_{\Ba}(x) \mid \Ba \in \CP_{n,e}\}$ gives a
$\ZZ$-basis of $\Xi^n$.  Put $\Xi_{\CC} = \CC\otimes \Xi$.  Then
$\{ p_{\Ba}(x) \mid \Ba \in \CP_{n,e}\}$ gives rise to a basis of
$\Xi_{\CC}$.
\para{1.3.}
Let $q, t$ be independent indeterminates and let $F = \CC(q,t)$
be the field of rational functions in $q,t$.  
We consider the $F$-algebra of symmetric functions 
$\Xi_F = F\otimes_{\ZZ}\Xi$ with coefficients in $F$.
Let $\Ba = (\a^{(0)}, \dots, \a^{(e-1)})$ be an $e$-partition.
For each partition 
$\a^{(k)} : \a^{(k)}_1 \ge \a^{(k)}_2 \ge \cdots $, put
\begin{equation*}
z_{\a^{(k)}}(q,t) = \prod_{j=1}^{l(\a^{(k)})}\frac {1-\z^kq^{\a_j^{(k)}} }
                                       {1 - \z^kt^{\a_j^{(k)}}}.
\end{equation*} 
We then define $z_{\Ba}(q,t) \in F$ by 
\begin{equation*}
\tag{1.3.1}
z_{\Ba}(q,t) = z_{\Ba}\prod_{k=0}^{e-1}z_{\Ba^{(k)}}(q,t),
\end{equation*}
where $z_{\Ba}$ is the order of the centralizer of $w_{\Ba}$
in $W$ (a representative of the conjugacy class of $W$ corresponding to 
$\Ba \in \CP_{n,e}$).  Explicitly, $z_{\Ba}$ is given as follows.  
For $\Ba \in \CP_{n,e}$, 
put $l(\Ba) = \sum_{k = 0}^{e-1}l(\a^{(k)})$.  For a partition 
$\a = (1^{n_1}, 2^{n_2}, \dots)$, put  
$z_{\a} = \prod_{i\ge 1}i^{n_i}n_i!$.  Then 
$z_{\Ba} = e^{l(\Ba)}\prod_{k=0}^{e-1}z_{\a^{(k)}}$.
\par
We define a sesquilinear form on $\Xi_F$ by 
\begin{equation*}
\tag{1.3.2}
\lp p_{\Ba},p_{\Bb}\rp = \d_{\Ba,\Bb}z_{\Ba}(q,t)
\end{equation*}
for $\Ba,\Bb \in \CP_{n,e}$.
Let $x = (x^{(k)}_j)$, $y = (y^{(k)}_j)$ ($0 \le k \le e-1)$ 
be two sequences of
infinitely many variables, and we define an infinite product
of $x$ and $y$ by
\begin{align*}
\tag{1.3.3}
\Pi(x,y;q,t) &= 
   \prod_{k=0}^{e-1}\prod_{i,j = 1}^{\infty}\prod_{r= 0}^{\infty}
    \frac{1 - tx_i^{(k-r-1)}y_j^{(k)}q^r}
         {1 - x_i^{(k-r)}y_j^{(k)}q^r}  \\
             &= 
    \prod_{k=0}^{e-1}\prod_{i,j = 1}^{\infty}\prod_{r= 0}^{\infty}
    \frac{1 - tx_i^{(k)}y_j^{(k+r+1)}q^r}
         {1 - x_i^{(k)}y_j^{(k+r)}q^r}.
\end{align*}   
(Here the upper indices of the variables in the formula should
be read $\mod e$).
Note that in the case where $q = 0$, the product $\Pi(x,y;q,t)$ 
reduces to the product $\Om(x,y;t)$ introduced in [S1, 2.5],
(or rather [S2, (5.7.1)], see the remark there).
\begin{lem}
$\Pi(x,y;q,t)$ has the following expansion.
\begin{equation*}
\Pi(x,y;q,t) = \sum_{\Ba}z_{\Ba}(q,t)\iv p_{\Ba}(x)\ol{p_{\Ba}(y)}
\end{equation*} 
where $\Ba$ runs over all the $e$-partitions of any size.
\end{lem}
\begin{proof}
Taking the $\log$ on both sides of the first formula of (1.3.3), 
\begin{equation*}
\begin{split}
\log\Pi(&x,y;q,t) \\ 
&= \sum_k\sum_{i,j}\sum_{r =0}^{\infty}\sum_{m=1}^{\infty}
\biggl\{ \frac 1 m (x_i^{(k-r)}y_j^{(k)}q^r)^m - 
   \frac {t^m} m (x_i^{(k-r-1)}y_j^{(k)}q^r)^m \biggr\}.
\end{split}
\end{equation*}
By making use of the equation 
\begin{equation*}
\frac 1 e\sum_{a =0}^{e-1}(\z^{k}\z^{-k'})^a = \d_{k,k'},
\end{equation*}
the above formula can be written as 
\begin{equation*}
\begin{split}
\log \Pi(x,y;q,t) = 
   \sum_{a=0}^{e-1}&\sum_{m \ge 1}\sum_{r \ge 0}\sum_{k,k',i,j} 
  \biggl\{\frac 1 {em} \z^{ar}\z^{a(k-r)}\z^{-ak'}(x_i^{(k-r)}y_j^{(k')}q^r)^m
 \\      
&- \frac {t^m}{em}\z^{a(r+1)}\z^{a(k-r-1)}\z^{-ak'}
(x_i^{(k-r-1)}y_j^{(k')}q^r)^m\biggr\}. 
\end{split}
\end{equation*}
It follows that
\begin{align*}
\log\Pi(x,y;q,t) &= \sum_{a=0}^{e-1}\sum_{m\ge 1}
                     \sum_{r \ge 0}\sum_{k,k',i,j}
    \frac{1 - \z^{a}t^m}{em}\z^{ar}q^{rm}(x_i^{(k)}y_j^{(k')})^m
                      \z^{ak}\z^{-ak'} \\
                 &= \sum_{a=0}^{e-1}\sum_{m \ge 1}\sum_{k,k',i,j}
    \frac 1 {em}\cdot\frac{1-\z^at^m}{1-\z^aq^m}
                    \z^{ak}\z^{-ak'}(x_i^{(k)}y_j^{(k')})^m  \\
                 &= \sum_{a=0}^{e-1}\sum_{m=1}^{\infty}
     \frac 1 {em}\cdot\frac{1-\z^at^m}{1-\z^aq^m}
                   p^{(a)}_m(x)\ol{p^{(a)}_m(y)}.
\end{align*}
Hence we have
\begin{align*}
\Pi(x,y;q,t) &= \prod_{a=0}^{e-1}\prod_{m=1}^{\infty}
          \exp\biggl\{\frac 1 {em}\cdot\frac{1-\z^at^m}{1-\z^aq^m}
               p^{(a)}_m(x)\ol{p^{(a)}_m(y)}\biggr\}  \\
             &= \sum_{\Ba}z_{\Ba}(q,t)\iv p_{\Ba}(x)\ol{p_{\Ba}(y)}.
\end{align*}
\end{proof}
\para{1.5.}
For a fixed $y^{(k)}_j$, we define a function $g^{(k)}_{m,\pm}(x;q,t)$
as the coefficient of $(y_j^{(k)})^m$ in  
\begin{equation*}
\tag{1.5.1}
\prod_i\prod_{r \ge 0}
        \frac{ 1 - tx_i^{(k\mp r \mp 1)}y_j^{(k)}q^r}
             { 1 - x_i^{(k\mp r)}y_j^{(k)}q^r}
        = \sum_{m \ge 0}g^{(k)}_{m,\pm}(x;q,t)(y_j^{(k)})^m
\end{equation*}
and put, for each $\Ba = (\a^{(k)}_j) \in \CP_{n,e}$, 
\begin{equation*}
g_{\Ba,\pm}(x;q,t) = \prod_{j, k}g^{(k)}_{\a_j^{(k)}, \pm}.
\end{equation*}
Then we have 
\begin{equation*}
\tag{1.5.2}
\Pi(x,y;q,t) = \sum_{\Ba}g_{\Ba,+}(x;q,t)m_{\Ba}(y)
             = \sum_{\Ba}m_{\Ba}(x)g_{\Ba,-}(y;q,t). 
\end{equation*}
In fact, by comparing the first formula of (1.3.3) and (1.5.1), we have 
\begin{align*}
\Pi(x,y;q,t) &= \prod_{k,j}\sum_{\a_j^{(k)} \ge 0}
      g^{(k)}_{\a_j^{(k)},+}(x;q,t)(y_j^{(k)})^{\a_j^{(k)}}  \\
             &= \sum_{\Ba}g_{\Ba,+}(x;q,t)m_{\Ba}(y).
\end{align*}
This shows the first equality.  If we compare the second formula
of (1.3.3) and (1.5.1) by replacing $x$ and $y$, the second 
equality is obtained in a similar way.
\par
Now by using a similar argument as in [M1, VI, 2.7], we see
that
\begin{equation*}
\tag{1.5.3}
\lp g_{\Ba,+}(x;q,t), m_{\Bb}(x)\rp = 
\lp m_{\Ba}(x), g_{\Bb,-}(x;q,t)\rp = \d_{\Ba,\Bb}.
\end{equation*}
In particular, the functions $g_{\Ba, \pm}(x;q,t)$ form a basis of $\Xi_F$ 
dual to $m_{\Ba}$. Hence $\{ g_n(x;q,t) \mid n \ge 0\}$ are
algebraically independent over $F$, and $\Xi_F = F[g_1, g_2, \dots]$.
\par
The following lemma can be proved in a similar way as [M1, VI, 2.13]
by using (1.5.2) and (1.5.3).
\begin{lem} 
Let $E^{\pm}: \Xi_F \to \Xi_F$ be two $F$-linear operators.  Then 
the following conditions are equivalent.
\begin{enumerate}
\item 
$\lp E^+f , g\rp = \lp f, E^-g\rp$ for any $f,g \in \Xi_F$ 
\item
$E^+_x\Pi(x,y;q,t) = E^-_y\Pi(x,y;q,t)$, where the suffix $x$
indicates the action of $E^{\pm}$ on the $x$ variables, and similarly 
for $y$.
\end{enumerate}
\end{lem} 
\para{1.7.}  
We shall give an explicit form of the function
$g_{\Ba,\pm}(x;q,t)$.
For this, we prepare some notation.
For $\mu = (\mu_0, \dots, \mu_{e-1}) \in \ZZ^e_{\ge 0}$, put
\begin{equation*}
f^{(k,i)}_{\mu,\pm}(x;q,t) = 
       \prod_{a=0}^{e-1}\prod_{j=1}^{\mu_a}
         \frac{x_i^{(k\mp a)} - tx_i^{(k\mp a \mp 1)}q^{e(j-1)}}
              {1 - q^{ej}}.
\end{equation*} 
Let $\CM_m$ be the set of sequences $\Bmu = (\mu^{(1)}, \mu^{(2)}, \dots)$ such 
that $\mu^{(i)} \in \ZZ^e_{\ge 0}$ and that $\sum_i |\mu^{(i)}| = m$.
For $\m = (\m_1, \m_2, \dots, \m_e) \in \ZZ^e_{\ge 0}$, put 
$n(\mu) = \sum_j (j-1)\mu_j$.  Then 
\begin{prop}
For each $m \ge 0$, we have
\begin{equation*}
\tag{1.8.1}
g^{(k)}_{m,\pm}(x;q,t) = \sum_{\Bmu \in \CM_m}\prod_i
                          f_{\mu^{(i)},\pm}^{(k,i)}(x;q,t)q^{n(\mu^{(i)})}.
\end{equation*}
\end{prop}  
\begin{proof}
By [M1, I, \S 2, Ex. 5], the following identity of formal power series
is known.
\begin{equation*}
\prod_{i=0}^{\infty}\frac{ 1 - bq^it}{1 - aq^it} = 
                 \sum_{m \ge 0}\biggl(\prod_{i=1}^m
                      \frac{a -bq^{i-1}}{1 - q^i\phantom{**}}\biggr)t^m.
\end{equation*}
Substituting this into 
\begin{align*}
A = 
\prod_{r \ge 0}\frac{1 - tx_i^{(k\mp r \mp 1)}yq^r}
                    {1 - x_i^{(k\mp r)}yq^r}
= \prod_{a=0}^{e-1}\prod_{r \ge 0}\frac{ 1-tx_i^{(k\mp a \mp 1)}
                                           yq^{re+a}}
                                       {1 - x_i^{(k\mp a)}yq^{re+a}}.
\end{align*}
with 
$a = x_i^{(k\mp a)}, b = tx_i^{(k\mp a\mp 1)}$, $t = q^ay$, 
$q = q^e$, we see
that 
\begin{align*}
A &= 
\prod_{a=0}^{e-1}\sum_{\mu_a =0}^{\infty}\prod_{j=1}^{\mu_a}
              \frac{x_i^{(k\mp a)} - tx_i^{(k\mp a \mp 1)}q^{e(j-1)}}
                   {1 - q^{ej}} (q^ay)^{\mu_a} \\
  &= \sum_{\mu \in \ZZ_{\ge 0}^e}f^{(k,i)}_{\mu, \pm}(x;q,t)q^{n(\mu)}y^{|\mu|}.
\end{align*}
It follows that
\begin{equation*}
\prod_{i \ge 1}\prod_{r \ge 0}
        \frac{ 1 - tx_i^{(k\mp r \mp 1)}y_j^{(k)}q^r}
             { 1 - x_i^{(k\mp r)}y_j^{(k)}q^r}
= \sum_{m \ge 0}\sum_{\Bmu \in \CM_m}\prod_{i \ge 1}
          f_{\mu^{(i)},\pm}^{(k,i)}(x;q,t)q^{n(\mu^{(i)})}(y_j^{(k)})^m.
\end{equation*} 
By comparing this with (1.5.1), we obtain the required formula.
\end{proof}
\par\medskip
{\bf Remark 1.9.}
\addtocounter{thm}{1}
$g_{m,\pm}^{(k)}(x;0,t)$ coincides with the function 
$q_{m,\pm}^{(k)}(x;t)$ introduced in [S1, 2.2].
By using a similar, but much simpler arguments as above, one 
obtains an alternative expression of $q_{m,\pm}^{(k)}(x;t)$
as follows.
\begin{equation*}
q_{m,\pm}^{(k)}(x;t) = \sum_{\mu \in \CP_{m}}
                |\FS_{\mu}|\iv\sum_{w \in \FS_{\mu}}
                  w\biggl\{ 
      \prod_{i=1}^{l(\mu)}(x_i^{(k)} - tx_i^{(k\mp 1)})
                  (x_i^{(k)})^{\mu_i-1}\biggr\}, 
\end{equation*}
where $\FS_{\mu}$ is the stabilizer of $\mu$ in $\FS_{m}$.
(Here we are considering finite variables 
$x_i^{(k)}, x_i^{(k\mp 1)}\,  (1 \le i \le m)$, and $\FS_m$ acts 
on both variables).
\para{1.10.}
The notion of symbols was introduced in [S1].  (Although a more
general setting was discussed in [S2], we do not use it 
in the discussion below.  We remark that similar symbols were
also considered by G. Malle in [Ma]).
Let $\Bm = (m_0, \dots, m_{e-1})$ be as before.  We denote by  
$Z_n^{0,0} = Z_n^{0,0}(\Bm)$ the set of $e$-partitions 
$\Ba = (\a^{(0)}, \dots, \a^{e-1)}) \in \CP_{n,e}$ 
such that each $\a^{(k)}$ is written ( as an element in $\ZZ^{m_k}$)
 in  the form $\a^{(k)}: \a_1^{(k)} \ge  \cdots \ge \a_{m_k}^{(k)} \ge 0$.
We express $\Ba$ as $\Ba = (\a^{(k)}_j)$ in matrix form.
Let us fix integers $r \ge s \ge 0$ and define  
an $e$-partition 
$\BLa^0 = \BLa^0(\Bm, s, r) = (\vL^{(0)}, \dots, \vL^{(e-1)})$ as follows.
\begin{align*}
\tag{1.10.1}
\vL^{(0)} &: (m_0-1)r \ge \cdots 2r \ge r \ge 0, \\
\vL^{(i)} &: s + (m_i-1)r \ge \cdots \ge s + 2r \ge s + r \ge s
\end{align*}
for $i = 0, \dots, e-1$.  We denote by 
$Z_n^{r,s} = Z_n^{r,s}(\Bm)$ the set of $e$-partitions of the form
$\BLa = \Ba + \BLa^0$, where $\Ba \in Z_n^{0,0}$ and the sum is taken 
entry-wise.  We write $\BLa = \BLa(\Ba)$ if $\BLa = \Ba + \BLa^0$, 
and call it the $e$-symbol of type $(r,s)$ corresponding to $\Ba$. 
We often denote the symbol $\BLa = (\vL^{(0)}, \dots, \vL^{(e-1)})$
in the form $\BLa = (\vL^{(k)}_j)$ with 
$\vL^{(k)}: \vL^{(k)}_1 > \dots > \vL^{(k)}_{m_k}$ 
for $k = 0, \dots, e-1$.
\par
Put $\Bm' = (m_0+1, \dots, m_{e-1} +1)$, and define 
a shift operation $Z_n^{r,s}(\Bm) \to Z_n^{r,s}(\Bm')$ by
associating 
$\BLa' = (\vL_0', \dots, \vL_{e-1}') \in Z_n^{r,s}(\Bm)$
to $\BLa = (\vL_0, \dots, \vL_{e-1})\in Z_n^{r,s}(\Bm)$,
where $\vL_0' = (\vL_0 + r) \cup \{ 0\}$, and  
$\vL_k' = (\vL_k + r) \cup \{s\}$ for $k = 0, \dots, e-1$.
In other words, for $\BLa = \BLa(\Ba)$, $\BLa'$ is obtained as
$\BLa' = \Ba + \BLa^0(\Bm',s, r)$,  where $\Ba$ is 
regarded as an element of $Z_n^{0,0}(\Bm')$ by adding 0 in the entries
of $\Ba$.
We denote by $\bar Z_n^{r,\Bs}$ the set of classes in 
$\coprod_{\Bm'}Z_n^{r,\Bs}(\Bm')$ under the equivalence relation generated
by shift operations.  
Note that $\CP_{n,e}$ coincides with the set $\bar Z_n^{0,0}$.  
Also note that $\BLa^0$ is regarded as a symbol in $Z_n^{r,s}$ with 
$n = 0$.
\par
Two elements $\BLa$ and $\BLa'$ in $\bar Z_n^{r,\Bs}$ are said 
to be similar, and are written as $\BLa \sim \BLa'$, 
if there exist representatives in $Z_n^{r,\Bs}(\Bm)$  such that
all the entries of them coincide with multiplicities. 
The set of symbols which are similar to a fixed symbol is called a
family in $Z_n^{r,\Bs}$.
\par
We define a function $a : Z_n^{r,\Bs} \to \ZZ_{\ge 0}$, for 
$\BLa \in  Z_n^{r,\Bs}$, by 
\begin{equation*}
\tag{1.10.2}
a(\BLa) = \sum_{\la,\la' \in \BLa}\min\{\la,\la'\}
            - \sum_{\m,\m' \in \BLa^0}\min\{\m,\m'\}.
\end{equation*}
The function $a$ on 
$Z_n^{r,\Bs}$ is invariant under the shift operation, and
it induces a function $a$ on $\bar Z_n^{r,\Bs}$. 
Clearly, the $a$-function takes a constant value on each family 
in $Z_n^{r,\Bs}$.
We regard the $a$-function as a function on $Z_n^{0,0}$ by using 
the bijection $Z_n^{0,0} \simeq Z_n^{r,s}$. 
\para{1.11.}
Hall-Littlewood functions $P^{\pm}_{\BLa}(x;t)$ and
$Q^{\pm}_{\BLa}(x;t)$ attached to 
symbols $\BLa$ were introduced in [S1].
We shall now construct a two parameter version of Hall-Littlewood 
functions. 
Let us introduce a total order $\Ba \prec \Bb$ on $Z_n^{0,0}$ 
such that $a(\Ba) \ge a(\Bb)$ whenever $\Ba \prec \Bb$ and that
each family in  $Z_n^{0,0}$ forms an interval.
\par
The following proposition is easily obtained by a similar 
argument as in Remark 4.9 in [S1] (i.e., a generalization of 
Gram-Schmidt orthogonalization process) if one notices that
$\Pi(x,y;0,t)$ coincides with $\Om(x,y;t)$ in 
[S1, 2.5]. 
\begin{prop} 
There exists a unique function $P^{\pm}_{\BLa}(x;q,t) \in \Xi_F$ 
for $\BLa \in Z_n^{r,s}$ satisfying the following two properties.
\begin{enumerate}
\item $P^{\pm}_{\BLa}(x;q,t)$ for $\BLa = \BLa(\Ba)$ can be 
expressed in terms of $s_{\Bb}(x)$ as 
\begin{equation*}
P^{\pm}_{\BLa} = s_{\Ba} + \sum_{\Bb}u^{\pm}_{\Ba,\Bb}s_{\Bb},
\end{equation*}
with $u_{\Ba,\Bb}^{\pm} \in F$, where $u^{\pm}_{\Ba,\Bb} = 0$
unless $\Bb \prec \Ba$ and $\Bb \not\sim \Ba$.
\item
$\lp P^+_{\BLa}, P^-_{\BLa'}\rp = 0$ unless $\BLa \sim \BLa'$.
\end{enumerate}
\end{prop}
\par
We then define $Q^{\pm}_{\BLa}(x;q,t)$ as the dual of 
$P^{\mp}_{\BLa}(x;q,t)$, i.e., by the property that
\begin{equation*}
\lp P^{+}_{\BLa}, Q^-_{\BLa'}\rp = 
       \lp Q^{+}_{\BLa}, P^-_{\BLa'}\rp = \d_{\BLa,\BLa'}.
\end{equation*}
$P^{\pm}_{\BLa}(x;q,t)$, $Q^{\pm}_{\BLa}(x;q,t)$ are called
the Macdonald functions associated to complex reflection groups 
$W$ (with respect to symbols in $Z_n^{r,s}$).
\par\medskip\noindent
{\bf Remark 1.13.}\ (i) The orthogonality relations of
Macdonald functions given above imply, by [M1, VI, 2.7], that
\begin{align*}
\tag{1.13.1}
\Pi(x,y;q,t) &= \sum_{\BLa}P^+_{\BLa}(x;q,t)Q^-_{\BLa}(y;q,t)  \\
             &= \sum_{\BLa}Q^+_{\BLa}(x;q,t)P^-_{\BLa}(y;q,t).
\end{align*} 
\par
(ii) In the case where $q = 0$, the scalar 
product given in (1.3.2) coincides with the one given in [S1, 4.7].
Then by Proposition 4.8 in [S1], one sees that
\begin{equation*}
\tag{1.13.2}
P^{\pm}_{\BLa}(x;0,t) = P^{\pm}_{\BLa}(x,t), \quad
Q^{\pm}_{\BLa}(x;0,t) = Q^{\pm}_{\BLa}(x;t),
\end{equation*} 
where the right hand sides are the Hall-Littlewood functions 
defined in [S1]. 
\par \ (iii) In the case where $q = t$, the scalar product in 
(1.3.2) coincides 
with the usual scalar product on the space $\Xi_{\QQ}$, where 
the Schur functions form an orthonormal basis of it.  Hence 
Proposition 1.12 implies that 
\begin{equation*}
\tag{1.13.3}
P^{\pm}_{\BLa(\Ba)}(x;t,t) = Q^{\pm}_{\BLa(\Ba)}(x;t,t) = s_{\Ba}(x).
\end{equation*}  
\vspace{1cm}
\section{Macdonald operators}
\para{2.1.}
The original Macdonald functions related to symmetric groups 
are characterized as the simultaneous eigenfunctions of 
Macdonald operators (see [M1, VI]).   
In this section, we shall construct certain operators which can be
viewed as a generalization of Macdonald operators to the case of
$W$.  Here we restrict ourselves to the case where symbols are
of the type $(r, s)$ with $r = 1$ and $s = 0$.
(We note that the arguments in this section can not be applied to 
other types of symbols. See Remark 2.9.)
In particular, $\BLa^0 = \Bdel$, where 
$\Bdel = (\d^{(0)}, \dots, \d^{(e-1)})$ with 
$\d^{(k)} = (m_k-1, \dots, 1,0)$. 
Hence, in the case where $e = 2$ (i.e., $W$ is the Weyl group 
of type $B_n$) with $m_1 = m_0+1$, 
these symbols are exactly the ones used to parameterize 
unipotent characters of 
finite classical groups $Sp_{2n}(\FF_q)$ or $SO_{2n+1}(\FF_q)$ by 
Lusztig [L1].  
Each family $\CF$ in $Z_n^{1,0}$ contains a unique element
$\BLa_{\CF} = (\vL^{(k)}_j)$ with the property 
\begin{equation*}
\tag{2.1.1}
\vL^{(0)}_1 \ge \vL^{(1)}_1 \ge \cdots 
       \ge \vL^{(e-1)}_1 \ge \vL^{(0)}_2 \ge \vL^{(1)}_2 \ge \cdots. 
\end{equation*}
Such an element is called a special symbol associated to the family 
$\CF$. The set of families is in bijection with the set
of special symbols.  Special symbols are regarded as partitions
of $N = \sum \vL_j^{(k)}$ by (2.1.1).
\par
We shall define a partial order on the set of families in 
$Z_n^{1,0}$.  Let $\CF$ and $\CF'$ be families in $Z_n^{1,0}$
and $\BLa_{\CF}, \BLa_{\CF'}$ be special symbols corresponding to them.
We put $\CF < \CF'$ if $\BLa_{\CF} < \BLa_{\CF'}$ with respect to the
dominance order on $\CP_N$.
Recall that for $\la = (\la_i), \mu = (\mu_i) \in \CP_N$, 
the dominance order $\la < \mu$ is defined by the condition that
\begin{equation*}
\sum_{i=1}^k\la_i \le \sum_{i=1}^k\mu_i, \qquad (1 \le k \le \sum_j m_j).
\end{equation*} 
We define a partial order on $Z_n^{1,0}$ by inheriting the partial
order on the set of families.
For $\la= (\la_i) \in \CP_N$, put 
\begin{equation*}
n(\la) = \sum_{i \ge 1}(i-1)\la_i.
\end{equation*}
Then, for each $\BLa$ in a family $\CF$, the value $a(\BLa)$ is given as
$\a(\BLa) = n(\BLa_{\CF}) - n(\BLa^0)$, where $\BLa_{\CF}$ is 
regarded as an element in $\CP_M$.
In particular, we have $a(\BLa) > a(\BLa')$ if $\BLa < \BLa'$.
As before, by using the bijection $Z_n^{0,0} \simeq Z_n^{1,0}$, 
we consider the partial order on $Z_n^{0,0}$, which will be denoted 
by the same symbol.
\para{2.2.}
In order to construct Macdonald operators, we shall start with
finitely many variables $x = x_{\Bm}$.
We consider the expansion of $\Pi(x,y;q,t)$ in the case of finitely
many variables.  Assume that $x = (x_j^{(k)}) = x_{\Bm}$ is as in 
1.2 and put  $\Xi_{\Bm, F} = F\otimes_{\ZZ}\Xi_{\Bm}$. 
Let us denote by $\CP_{\Bm}$ the set of $e$-partitions  
$\Ba = (\a_j^{(k)})$ such that $l(\a^{(k)}) \le m_k$. 
Then $m_{\Ba}(x) = 0$ unless $\Ba \in \CP_{\Bm}$, and those non-zero 
$m_{\Ba}(x)$ form a basis of $\Xi_{\Bm,F}$.  The same is true for
Schur functions.  Also by 
a similar argument as in [M1] we see that $g_{\Ba}(x;q,t)$   
such that $\Ba \in \CP_{\Bm}$ form an $F$-basis of $\Xi_{\Bm,F}$.
\par
Now by substituting $x^{(k)}_j = y^{(k)}_j = 0$ for $j > m_k$, we have
a finite version of (1.5.2), i.e., for $x = x_{\Bm}, y = y_{\Bm}$, we have
\begin{equation*}
\tag{2.2.1}
\Pi(x,y;q,t) = \sum_{\Ba \in \CP_{\Bm}}g_{\Ba,+}(x;q,t)m_{\Ba}(y)
             = \sum_{\Ba \in \CP_{\Bm}}m_{\Ba}(x)g_{\Ba,-}(y;q,t). 
\end{equation*}
This enables us to define a scalar product on $\Xi_{\Bm,F}$ by 
\begin{equation*}
\tag{2.2.2}
\lp g_{\Ba,+}, m_{\Bb}\rp = \lp m_{\Ba}, g_{\Bb,-}\rp = \d_{\Ba,\Bb}.
\end{equation*}
\par
We now consider the restriction of the functions 
$P^{\pm}_{\BLa}$ to $\Xi_{\Bm,F}$. By Proposition 1.12, one sees that 
$\{ P^{\pm}_{\BLa(\Ba)} \mid \Ba \in \CP_{\Bm}\}$ form a basis 
of $\Xi_{\Bm,F}$.  Moreover, the finite variables version of
Proposition 1.12 holds, and $P^{\pm}_{\BLa} \in \Xi_F$ is obtained
as the limit of $P^{\pm}_{\BLa} \in \Xi_{\Bm,F}$ determined by these 
properties. 
\para{2.3.}
Let 
\begin{equation*}
\CI = \CI(\Bm) = \{ \Bi = \begin{pmatrix}
                   i_0  \\
                   \vdots \\
                   i_{e-1}
                          \end{pmatrix}  \in \ZZ^e 
              \mid 1 \le i_k \le m_k \}.
\end{equation*}
For each $\Bi \in \CI$ and $u \in F$,  
we define an $F$-linear operator $T^{\pm}_{u, \Bi}: F[x] \to F[x]$ by 
$T^{\pm}_{u,\Bi}f = f'$, where $f'$ is a polynomial obtained from $f$ by 
replacing the variables $x_{i_k}^{(k)}$ by $ux_{i_{k\mp 1}}^{(k\mp 1)}$
for $k = 0, \dots, e-1$.
(Here we understand that $i_e = i_o$).
More generally, for each $r$ such that $1 \le r \le M_1 = \min_k \{m_k\}$, 
we define $\CI_r$ as the set of 
$J = \{\Bi_1, \dots, \Bi_r\}$ consisting of $\Bi_k \in \CI$ such that
any two $\Bi_{k}$ have no common entries (i.e., $\Bi_j - \Bi_k$ does
not contain 0 entries for each $j \ne k$ as vectors in $\ZZ^e$).
Then we define, for each $J \in \CI_r$, an operator 
$T^{\pm}_{u,J}: F[x] \to F[x]$ by 
$T^{\pm}_{u,J} = \prod_{k =1}^rT^{\pm}_{u,\Bi_k}$.
Note that $T^{\pm}_{u,\Bi_k}$ in the product commute with each other, and so 
$T^{\pm}_{u,J}$ does not depend on the order of the product.  
\par
Let $Z = Z(\Bm)$ be the set of sequences $\Bb = (\b^{(k)}_j)$ 
($0 \le k \le e-1, 1 \le j \le m_k$) with 
$\b^{(k)}_i \in \ZZ_{\ge 0}$.  
For each $\Bb \in Z$, we denote by $[\Bb]$ the element in 
$Z$ obtained from $\Bb$ by permuting the entries inside each row, 
so that each row
is arranged in decreasing order.
We often regard $\Bb$ as a matrix, 
and denote its $k$-th row by $\b^{(k)}$. 
Then the set $J$ is regarded as a subset of the set of indices 
$\{ (k,j) \mid 0 \le k \le e-1, 1 \le j \le m_k\}$ of $(\b^{(k)}_j)$
\par
We put,
for each $\Bb \in Z$ and $J \in \CI_r$, 
\begin{equation*}
\lp \Bb, J\rp = \sum_{(k,j) \in J}\b_{j}^{(k)}.
\end{equation*}
$\FS_{\Bm}$ acts naturally on $Z$ and on $\CI_r$, respectively, 
and this pairing is $\FS_{\Bm}$-invariant.  
We also note that the action of $\FS_{\Bm}$ on $\CI_r$ is transitive. 
\par
The operation of $T^{\pm}_{q,J}$ on $x = (x_i^{(k)})$ also induces 
an action $\Bb \mapsto \Bb_{J\pm}$ on $Z$ by permuting the entries of $\Bb$
so that $T^{\pm}_{q,J}x^{\Bb} = q^{\lp \Bb,J\rp}x^{\Bb_{J\pm }}$,
where, as usual, $x^{\Bb}$ denotes the monomial 
$\prod_{i,k}(x_i^{(k)})^{\b_i^{(k)}}$.
For each $\Bb \in Z$, we define a function $a_{\Bb}(x)$ by 
\begin{equation*}
a_{\Bb}(x) = \sum_{w \in \FS_{\Bm}}\ve(w)w(x^{\Bb}) 
\end{equation*}
\par
We now define, for each $1 \le r \le M_1$,  an $F$-linear
operator
$D^{r}_{\pm}(q,t)$ on $F[x]$ by
\begin{equation*}
\tag{2.3.1}
D^{r}_{\pm}(q,t)
= a_{\Bdel}(x)\iv \sum_{w \in \FS_{\Bm}}
     \ve(w)\sum_{J \in \CI_r}x^{w(\Bdel)_{J\pm}}
                 t^{\lp w(\Bdel), J\rp}T^{\pm}_{q,J}. 
\end{equation*}
Then $D^r_{\pm}(q,t)$ can be written as
\begin{equation*}
\tag{2.3.2}
D^r_{\pm}(q,t) = \sum_{J \in \CI_r}A_J^{\pm}(x;t)T_{q,J}^{\pm}, 
\end{equation*}
with
\begin{align*}
A_J^{\pm}(x;t) &= a_{\Bdel}(x)\iv\sum_{w \in \FS_{\Bm}}
                   \ve(w)x^{w(\Bdel)_{J\pm}}t^{\lp w(\Bdel),J\rp} \\
               &= a_{\Bdel}(x)\iv T^{\pm}_{t,J}a_{\Bdel}(x)  \\
               &= t^{r(r-1)/2}
     \prod_{k=0}^{e-1}\prod_{\substack{(k,i) \notin J \\
                                        (k,j) \in J }}
       \frac{\ x_i^{(k)} - tx_{j'}^{(k\mp 1)}}{x_i^{(k)} - x_j^{(k)}},
\tag{2.3.3}
\end{align*}
where we write $J = \{ \Bi_1, \dots, \Bi_r\}$, and 
take $(k\mp 1, j') \in \Bi_a$ if $(k,j) \in \Bi_a$ for
$a = 1, \dots, r$.
The formulas (2.3.2) and (2.3.3) are the analogue of  
$(3.4)_r$ and $(3.5)_r$ in [M, VI, 3].  
\par
First we show that 
\begin{lem}
\begin{enumerate}
\item
For $\Ba \in \CP_{\Bm}$, we have
\begin{equation*}
D^{r}_{\pm}(q,t)m_{\Ba}(x)  
     = \sum_{\Bb}
          \sum_{J \in \CI_r}t^{\lp\Bdel, J\rp}
                 q^{\lp \Bb, J\rp }
                     s_{(\Bb +\Bdel)_{J\pm} - \Bdel}(x), 
\end{equation*} 
where $\Bb \in Z$ runs over all the row permutations of $\Ba$.
\item
For $\Ba \in \CP_{\Bm}$, we have, 
\begin{equation*}
D^{r}_{\pm}(q,t)m_{\Ba}(x) = 
   \sum_{\Bb \in Z_n^{0,0}}b^{r,\pm}_{\Ba,\Bb}(q,t)m_{\Bb}(x)
\end{equation*}
with $b^{r,\pm}_{\Ba,\Bb}(q,t) \in F$, 
where $b^{r,\pm}_{\Ba,\Bb}(q,t) = 0$ unless $\Bb \sim \Ba$ or $\Bb < \Ba$.
In particular, $D^r_{\pm}$ is an operator on the space $\Xi_{\Bm,F}$.
\end{enumerate}
\end{lem}
\begin{proof}
For each $\Bb \in Z_n^{0,0}$, we have 
\begin{equation*}
D^{r}_{\pm}(q,t)x^{\Bb}  
= a_{\Bdel}(x)\iv\sum_{w_1 \in \FS_{\Bm}}
     \ve(w_1)\sum_{J \in \CI_r}t^{\lp w_1(\Bdel), J\rp}
           q^{\lp \Bb, J\rp}
               x^{(\Bb + w_1(\Bdel))_{J\pm }}.
\end{equation*}
If we replace $\Bb$ by $w_2(\Bb)$ for $w_2 \in \FS_{\Bm}$ and put  
$w_2 = w_1w$, the term $(\Bb + w_1(\Bdel))_{J\pm }$ 
(resp. $\lp\Bb, J\rp$ )
is replaced by $w_1((w(\Bb) + \Bdel)_{J'\pm})$ 
(resp. $\lp w(\Bb), J'\rp$ ), respectively, with $J' = w_1\iv(J)$. 
It follows that for $\Ba \in Z_n^{0,0}$, we have
\begin{align*}
D^{r}_{\pm}(q,t)m_{\Ba}(x) &= |\FS_{\Ba}|\iv a_{\Bdel}(x)\iv  
\sum_{w,w_1 \in \FS_{\Bm}}\ve(w_1) \\
    &\phantom{*******}\times\sum_{J' \in \CI_r}t^{\lp \Bdel, J'\rp}
           q^{ \lp w(\Ba), J'\rp}
              x^{w_1(( w(\Ba) + \Bdel)_{J'\pm})} \\
    &= |\FS_{\Ba}|\iv\sum_{w \in \FS_{\Bm}}
          \sum_{J \in \CI_r}t^{\lp\Bdel, J\rp}
                 q^{\lp w(\Ba), J\rp }
                     s_{(w(\Ba) +\Bdel)_{J\pm} - \Bdel}(x).
\end{align*}
This proves (i).
\par
Next we show (ii).  Since $\Bdel = \BLa^0$ by our assumption, 
$\Ba + \Bdel$ coincides with the symbol $\BLa(\Ba)$ for 
$\Ba \in Z_n^{0,0}$. We note that if $w \neq 1$, then the 
symbol $w(\Ba) + \Bdel$ belongs to a family strictly smaller than
the family containing $\Ba$.  On the other hand, if $w = 1$, 
 $(\Ba +\Bdel)_{J\pm}$ is obtained
from the symbol $\Ba +\Bdel$ by permuting some entries, 
and $s_{(\Ba+\Bdel)_{J\pm} -\Bdel}$ coincides with 
$\pm s_{\Bg}$, where $\Bg +\Bdel$ is obtained from 
$(\Ba+\Bdel)_{J\pm}$ by rearranging the rows in decreasing 
order. 
It follows that $\Bg + \Bdel$ is contained in the same family as 
$\Ba + \Bdel$. Hence, for $\Bb \in Z_n^{0,0}$ in the expression  of 
(ii), we see that 
$\Bb < \Ba$ if $w \ne 1$ and $\Bb \sim \Ba$ if $w = 1$.
\end{proof} 
\par
Next we show
\begin{lem}
The operators $D^{r}_{\pm}$ are adjoint each other, i.e., we have 
\begin{equation*}
\tag{2.5.1}
\lp D^r_+f, g\rp = \lp f, D^{r}_{-}g\rp\, ,  
                   \quad (f, g \in \Xi_{\Bm,F}).
\end{equation*}
\end{lem}
\begin{proof}
By Lemma 1.6, (2.5.1) is equivalent to the formula
\begin{equation*}
\tag{2.5.2}
\Pi\iv (D^{r}_{+})_x\Pi = \Pi\iv (D^{r}_{-})_y\Pi.
\end{equation*}
But for $J \in \CI_r$, we have 
\begin{align*}
\Pi\iv (T^{+}_{q,J})_x\Pi &= \prod_{k=0}^{e-1}\prod_{j\ge 1}
                               \prod_{(k,i) \in J}
   \frac{1 - x_{i}^{(k)}y_j^{(k)}}{1 - tx_{i}^{(k)}y_j^{(k + 1)}}, \\
\Pi\iv (T^{-}_{q,J})_y\Pi &= \prod_{k=0}^{e-1}\prod_{i\ge 1}
                                \prod_{(k,j) \in J}
   \frac{1 - x_{i}^{(k)}y_j^{(k)}}{1 - tx_{i}^{(k-1)}y_j^{(k)}}. \\
\end{align*}
It follows that both of $\Pi\iv (D^{r}_{+})_x\Pi$ and 
$\Pi\iv (D^{r}_{-})_y\Pi$ 
are independent of $q$.
Hence in the proof of (2.5.2), we may assume that $q = t$. 
In other words, we have only to prove (2.5.1) under the assumption 
that $q= t$.
\par
Now assume that $q = t$.  Since 
\begin{equation*}
T^{\pm}_{t,J}(x^{w(\Bdel)}f) = x^{w(\Bdel)_{J\pm}}
      t^{\lp w(\Bdel), J\rp}T^{\pm}_{t,J}f
\end{equation*}
for any polynomial $f \in F[x]$, we have
\begin{equation*}
D^{r}_{\pm}(t,t)f = a_{\Bdel}\iv\sum_{J \in \CI_r}
       T^{\pm}_{t,J}(a_{\Bdel}f).
\end{equation*}
It follows that for any $\Ba \in Z_n^{0,0}$, 
\begin{align*}
D^{r}_{\pm}(t,t)s_{\Ba} &= 
    a_{\Bdel}\iv\sum_{J \in \CI_r}T^{\pm}_{t,J}(a_{\Ba + \Bdel})  \\
 &= \sum_{J \in \CI_r}t^{\lp \Ba + \Bdel, J\rp}
        s_{(\Ba + \Bdel)_{J\pm}-\Bdel}.
\end{align*}
As before,  $s_{(\Ba+\Bdel)_{J\pm}- \Bdel}$ coincides with $\pm s_{\Bb}$, where 
$\Bb + \Bdel = [(\Ba+\Bdel)_{J\pm}]$ (under the notation in 2.3).    
Now in the case where $q=t$, $\{ s_{\Ba}(x) \mid \Ba \in \CP_{\Bm}\}$ 
is an orthonormal basis of $\Xi_{\Bm,F}$.  It follows that
\begin{equation*}
\tag{2.5.3}
\lp D^{r}_{+}(t,t)s_{\Ba}, s_{\Bb}\rp = 
       \sum_{J}\ve_{J+}t^{\,\lp \Ba+\Bdel, J\rp}, 
\end{equation*}
where $J$ runs over all the elements in $\CI_r$ such that
$[(\Ba+\Bdel)_{J+}]$ coincides with $\Bb+\Bdel$, and 
$\ve_J = (-1)^{\l(w)}$ with $w \in \FS_{\Bm}$ such that 
$[(\Ba+\Bdel)_{J+}] = w((\Ba+\Bdel)_{J+})$. 
But if $\Bb +\Bdel = w((\Ba +\Bdel)_{J+})$, we have
$\Ba + \Bdel = w\iv((\Bb + \Bdel)_{J'-}) = 
[(\Bb + \Bdel)_{J'-}]$  with $J' = w(J)$.   
Also in this case, 
\begin{equation*}
\lp \Bb+\Bdel, J'\rp = \lp\, w(\Ba+\Bdel)_{J'+}, J'\rp
                    = \lp w(\Ba+\Bdel), J'\rp 
                    = \lp \Ba+\Bdel, J\rp.
\end{equation*}
Since $\ve_{J+} = \ve_{J'-}$, we see that the right hand side
of (2.5.3) is equal to 
\begin{equation*}
\sum_{J'}\ve_{J'-}t^{\,\lp \Bb+\Bdel, J'\rp},
\end{equation*}
where $J' \in \CI_r$ runs over all the elements such that
$[(\Bb + \Bdel)_{J'-}] = \Ba +\Bdel$.
Clearly this coincides with $\lp s_{\Ba}, D^r_{-}(t,t)s_{\Bb}\rp$.
So the lemma is proved.
\end{proof}
\para{2.6.}
We fix a total order $\prec$ on $Z_n^{0,0}$ as in 1.11, so that it is 
compatible with the partial order $<$.  
Let $B^r_{\pm} = (b^{r,\pm}_{\Ba,\Bb})$ be the matrix consisting of 
the coefficients in the formula in Lemma 2.4  
(ii) with respect to the total order $\prec$.
We consider $B^r_{\pm}$ as a block matrix with respect to the equivalence
relation $\Ba \sim \Bb$, and denote it as $(B^{r,\pm}_{\CF,\CF'})$, where
$B^{r,\pm}_{\CF,\CF'}$ is the submatrix of $B^r_{\pm}$ corresponding 
to the families $\CF,\CF'$.  Then Lemma 2.4 (ii) implies that
$B^r_{\pm}$ is lower triangular as a block matrix.
We consider the Macdonald functions $P^{\pm}_{\BLa}$ constructed 
via $\prec$.  The following result shows that the set of 
Macdonald functions attached to symbols in a fixed family 
behaves as an eigenfunction for Macdonald operators, where the
eigenvalues should be replaced by the diagonal blocks 
$B^{r,\pm}_{\CF,\CF}$ of $B_{\pm}$.
\begin{prop} 
Let $P^{\pm}_{\BLa}(x;q,t) \in \Xi_{\Bm,F}$ be Macdonald functions attached 
to $\BLa = \BLa(\Ba)$.  Then we have
\begin{equation*}
D^{r}_{\pm}P^{\pm}_{\BLa} = \sum_{\Bb \sim \Ba}
                            b^{r,\pm}_{\Ba, \Bb}P^{\pm}_{\BLa(\Bb)},
\end{equation*}  
where the coefficients $b^{r,\pm}_{\Ba,\Bb}$ are the same as in Lemma 2.4 (ii).
\end{prop} 
\begin{proof}
By Proposition 1.12, $P^{\pm}_{\BLa}$ can be written as
\begin{equation*}
P^{\pm}_{\BLa} = m_{\Ba} + \sum_{\Bb \prec \Ba, \Bb \not\sim \Ba}
                    u'_{\Ba,\Bb}m_{\Bb}
\end{equation*}
with $u'_{\Ba,\Bb} \in F$.  It follows, by Lemma 2.4 (ii), that 
one can write as 
\begin{equation*}
D^r_{\pm}P^{\pm}_{\BLa}
   = \sum_{\Bb \sim \Ba}b^{r,\pm}_{\Ba, \Bb}P^{\pm}_{\BLa(\Bb)}
          + \sum_{\BLa' \prec \BLa, \BLa' \not\sim \BLa}
                 {c^{\pm}_{\BLa, \BLa'}}P^{\pm}_{\BLa'}.
\end{equation*}
Hence, for each $\BLa'$ such that $\BLa' \prec \BLa$ and that 
$\BLa' \not\sim \BLa$, we have  
\begin{equation*}
\lp D^r_+P^{+}_{\BLa}, P^{-}_{\BLa'}\rp 
     = c^+_{\BLa,\BLa'}.
\end{equation*}
On the other hand thanks to Lemma 2.5,  we have   
\begin{equation*}
\lp D^r_+P^+_{\BLa}, P^-_{\BLa'}\rp = 
         \lp P^+_{\BLa}, D^r_-P^-_{\BLa'}\rp  = 0
\end{equation*}
since $D^r_-P^-_{\BLa'}$ is a linear combination of 
$P^-_{\BLa''}$, where $\BLa'' \sim \BLa'$ or $\BLa'' \prec \BLa'$.
It follows that $c^+_{\BLa, \BLa'} = 0$ and the proposition holds 
for the $+$ case.  The $-$ case is similar. 
\end{proof}
\par
In view of Proposition 2.7, it is important to know the diagonal 
part of $B_{\pm}$.  By lemma 2.4, the matrix 
$B^{\pm}_{\CF,\CF}$ is described as follows.
\begin{lem} 
For $\Ba, \Bb \in Z_n^{0,0}$ such that $\Bb \sim \Ba$, we have
\begin{equation*}
b^{r,\pm}_{\Ba,\Bb}(q,t) = 
    \sum_{\substack{ J \in \CI_r \\
                        [\BLa_{J\pm}] = \BLa(\Bb)}}
             \ve_{\BLa,J\pm}(tq\iv)^{\lp\Bdel.J\rp}q^{\lp\BLa,J\rp},
\end{equation*}
where $\BLa = \Ba +\Bdel \in Z_n^{1,0}$ and 
$\ve_{\BLa, J\pm} = (-1)^{l(w)}$ for $w \in \FS_{\Bm}$
such that $[\BLa_{J\pm}] = w(\BLa_{J\pm})$.
\end{lem}
\addtocounter{thm}{1}
\noindent
{\bf Remark 2.9.}
The results in this section work only for a special type 
of symbols.  It seems to be difficult to extend the definition 
of Macdonald operators in (2.3.1) directly to a more general case.
A naive idea for the general situation is to replace 
$\Bdel$ by $\BLa^0$ in the definition (2.3.1).  
Then one gets some operator related to the symbols associated with 
$\BLa^0$.  However, the thus obtained operator does not
preserve the set of polynomials in $x$, since it involves the factor 
$a_{\BLa^0}(x)\iv$.  
If one leaves the denominator $a_{\Bdel}(x)$ unchanged, and replaces
$\Bdel$ in all other places, then the operator preserves 
polynomials, but it does not preserve the degree of them, and
is not so useful.
\par\medskip
\section{A characterization of Macdonald functions}
\para{3.1.}
We write the operator $D^{r}_{\pm}$ as $D^{r}_{\Bm,\pm}$ 
to indicate the dependence on $\Bm$. 
The operators $D^{r}_{\Bm,\pm}$ are not compatible with 
the restriction homomorphisms 
$\r_{\Bm',\Bm}$.
In the case where $r = 1$, one can
modify $D^1_{\Bm,\pm}$ as discussed in [M1], so that they are
compatible with $\r_{\Bm',\Bm}$.
Let us define an operator $E^{\pm}_{\Bm} = E^{\pm}_{\Bm}(q,t)$ 
on $\Xi_{\Bm, F}$ by
\begin{equation*}
\tag{3.1.1}
E^{\pm}_{\Bm} = t^{-M}D^1_{\Bm,\pm} - 
         \sum_{\Bi \in \CI}t^{\lp\Bdel,\Bi\rp - M}.              
\end{equation*} 
We show the following lemma.
\begin{lem} 
The operators $E^{\pm}_{\Bm}$ are compatible with $\r_{\Bm', \Bm}$, 
i.e., we have
\begin{equation*}
\r_{\Bm',\Bm}\circ E^{\pm}_{\Bm'} = E^{\pm}_{\Bm}\circ \r_{\Bm',\Bm}
\end{equation*}
for $\Bm' = (m_0+1, \dots, m_{e-1}+1)$.  
\end{lem}
\begin{proof}
Let us define another operator 
$\wt E^{\pm}_{\Bm} : \Xi_{\Bm,F} \to \Xi_{\Bm,F}$ by 
\begin{equation*}
\wt E^{\pm}_{\Bm} = t^{-M}D^1_{\Bm,\pm} - 
         \sum_{\Bi \in \CI}t^{\lp\Bdel,\Bi\rp - M}
              a\iv_{\Bdel}T^{\pm}_{1,\Bi}a_{\Bdel}.
\end{equation*} 
Then for each $\Ba \in \CP_{\Bm}$, we have
\begin{equation*}
\tag{3.2.1}
\wt E^{\pm}_{\Bm}m_{\Ba} = 
     \sum_{\Bb}\sum_{\Bi \in \CI}
          (q^{\lp\Bb,\Bi\rp}-1)t^{-\sum i_k}
               s_{(\Bb +\Bdel)_{\Bi\pm}-\Bdel},
\end{equation*}
where $\Bb \in Z$ runs over all the row permutations of $\Ba$, 
and $\Bi = (i_0, \dots, i_{e-1})$.
In fact, by applying [M1, VI, 4] one can write
$m_{\Ba} = \sum_{\Bb}s_{\Bb}$ with $\Bb$ as above.  Since 
$a\iv_{\Bdel}T_{1,\Bi}^{\pm}a_{\Bdel}s_{\Bb} = 
         s_{(\Bb +\Bdel)_{\Bi\pm}-\Bdel}$, 
we obtain (3.2.1).
\par
We claim that $\wt E^{\pm}_{\Bm}$ is compatible with 
$\r_{\Bm',\Bm}$, i.e., 
\begin{equation*}
\tag{3.2.2}
\r_{\Bm',\Bm}\circ  \wt E^{\pm}_{\Bm'} = \wt E^{\pm}_{\Bm}\circ \r_{\Bm',\Bm}.
\end{equation*}
Recall that the map 
$\r_{\Bm',\Bm} : \Xi_{\Bm',F} \to \Xi_{\Bm,F}$ is defined by 
substituting $x^{(k)}_{m_k+1} = 0$ for $k = 0, \dots, e-1$.
Take $\Ba \in Z_n^{0,0}(\Bm)$ and let $\Ba'\in Z_n^{0,0}(\Bm')$
be the element obtained by adding 0's to the last part
of $\Ba$.  
We consider the expression of $\wt E_{\Bm'}^{\pm}m_{\Ba'}$ as given in 
(3.2.1).
Let $\Bdel'$ be the element for $\Bm'$
corresponding to $\Bdel$ for $\Bm$, and take 
$\Bb' \in Z(\Bm')$.
We note that 
$s_{(\Bb' +\Bdel')_{\Bi\pm}-\Bdel'}$ goes to zero under
$\r_{\Bm',\Bm}$ if there exists 
some $k$ such that $\b^{(k)}_{m_k+1} \ne 0$ for $\Bb' = (\b^{(k)}_j)$. 
In fact, if we write $\Bb'+\Bdel' = (\g^{(k)}_j)$, then $\g^{(k)}_j >0$
for $j \le m_k$.  Hence if $\b^{(k)}_{m_k+1} \ne 0$ for some $k$,
then $(\Bb'+\Bdel')_{\Bi\pm} = ({\g'}^{(k)}_j)$ contains 
a row ${\g'}^{(k')}$ whose entries are all non-zero,   
and so $s_{(\Bb'+\Bdel')_{\Bi\pm}-\Bdel'}$ goes to zero.
\par
It follows that in the expression of $\wt E^{\pm}_{\Bm'}m_{\Ba'}$ 
in (3.2.1),
we may only consider $\Bb'$ (a row permutation of $\Ba'$ ) such that
the last column consists of zeros.  
One can embed $Z(\Bm)$ into $Z(\Bm')$ by adding 0 to the last part 
of $\Bb \in Z(\Bm)$.  Then those $\Bb'$ are identified, under
the embedding $Z_n^{0,0}(\Bm) \hookrightarrow Z_n^{0,0}(\Bm')$, with
a row permutation $\Bb$ of $\Ba$.
Now $\CI(\Bm)$ is also embedded into $\CI(\Bm')$ in the same way. 
Take $\Bi \in \CI(\Bm')$ such that $\Bi \not\in \CI(\Bm)$.
If $\Bi$ is not equal to  $\Bi_1 = {}^t(m_0+1, \dots, m_{e-1}+1)$, then 
$s_{(\Bb'+\Bdel')_{\Bi\pm}-\Bdel'}$ goes to zero under $\r_{\Bm',\Bm}$
by the same reason as before. 
If $\Bi = \Bi_1$, then $\lp \Bb', \Bi\rp = 0$.  In any case, such 
$\Bi$ does not give a contribution, and we may only consider 
$\Bi \in \CI(\Bm)$ in (3.2.1).  Hence (3.2.2) holds. 
\par
Let $\wt H_{\pm}$ be the coefficient matrix of $\wt E^{\pm}_{\Bm}m_{\Ba}$ in
terms of $m_{\Bb}$.  Then by Lemma 2.4 (ii), $\wt H_{\pm}$ can be
expressed as a block matrix $\wt H_{\pm} = (\wt H^{\pm}_{\CF,\CF'})$, where 
$\wt H^{\pm}_{\CF,\CF'} = 0$ unless $\CF = \CF'$ or $\CF' < \CF $.
Then by comparing (3.2.1) with Lemma 2.4 (i) on the diagonal parts, 
we see that
\begin{equation*}
\wt H^{\pm}_{\CF,\CF} = t^{-M}B^{\pm}_{\CF,\CF} - 
                             \sum_{\Bi\in \CI}t^{\lp\Bdel,\Bi\rp -M}. 
\end{equation*}
On the other hand, if we write 
$\wt H^{\pm}_{\CF,\CF} = (h^{\pm}_{\Ba,\Bb})_{\Ba,\Bb \in \CF}$, 
Proposition 2.7 together with (3.1.1) implies that
\begin{equation*}
E^{\pm}_{\Bm}P^{\pm}_{\BLa(\Ba)} = \sum_{\Bb \sim \Ba}
     h^{\pm}_{\Ba,\Bb}P^{\pm}_{\BLa(\Bb)}.
\end{equation*} 
But by (3.2.2), the matrix $\wt H^{\pm}_{\CF,\CF}$ does not depend 
on the shift operation under $\Bm' \to \Bm$.  This shows that
the operator $E^{\pm}_{\Bm}$ is compatible with $\r_{\Bm',\Bm}$.
The lemma is proved. 
\end{proof}
\para{3.3.}
By Lemma 3.2, one can define an operator $E^{\pm} = E^{\pm}(q,t)$ 
on $\Xi_F$ as the limit of $E^{\pm}_{\Bm}$.
Since $E^{\pm}_{\Bm}$ satisfies a similar formula as given
in Lemma 2.5, the operator $E^{\pm}$ also satisfies the adjointness
property, i.e., we have
\begin{equation*}
\tag{3.3.1}
\lp E^+f,g\rp = \lp f, E^-g\rp, \quad (f,g \in \Xi_F).
\end{equation*}  
By Lemma 2.4 (ii), one can write, for each $\Ba \in \CP_{n,e}$, 
\begin{equation*}
\tag{3.3.2}
E^{\pm}(q,t)m_{\Ba}(x) = \sum_{\Bb \in Z_n^{0,0}}
          h^{\pm}_{\Ba,\Bb}(q,t)m_{\Bb}(x)
\end{equation*}
with $h^{\pm}_{\Ba,\Bb}(q,t) \in F$, where 
$h^{\pm}_{\Ba,\Bb}(q,t) = 0$ unless $\Bb \sim \Ba$ or $\Bb < \Ba$.
Let $H_{\pm} = (h^{\pm}_{\Ba,\Bb})$, and write it as a block matrix
$H_{\pm} = (H^{\pm}_{\CF,\CF'})$ as in 2.6. 
\par
In the case of symmetric groups, the matrix $H_{\pm}$ is a triangular
matrix, with distinct eigenvalues.  This property was used to 
characterize Macdonald functions as eigenfunctions of Macdonald
operators.  As an analogy, it is likely that the following property 
holds for the diagonal parts of the block matrix $B_{\pm}$.
\par\medskip\noindent
{\bf Conjecture A.}\ Let $\CF, \CF'$ be any distinct families 
       in $Z_n^{0,0}$.  Then
       the matrices $H^{\pm}_{\CF,\CF}$ and $H^{\pm}_{\CF',\CF'}$ 
       have no common eigenvalues (according to the sign $+$ or $-$, 
       respectively). 
\par\medskip
We have verified the conjecture in the case where $e= 2$, and $n \le 5$.
\para{3.4.}
Before giving a characterization of Macdonald functions in terms 
of Macdonald operators, we prepare an easy lemma.
Let $A = (a_{ij})$ (resp. $B = (b_{ij})$) be a square matrix
of degree $m$ (resp. $n$), and 
let $C = (C_{\Ba,\Bb})_{1 \le \a,\b \le n}$ be a block matrix of
size $mn$, consisting of blocks $C_{\a,\b}$ of size $m$, defined by 
\begin{equation*}
C_{\Ba,\Bb} = \begin{cases}
        A - b_{\a,\a}I_{m} &\quad\text{ if } \b = \a, \\
        -b_{\a,\b}I_m       &\quad\text{ otherwise.}
              \end{cases}
\end{equation*}  
We consider a matrix equation 
$AX = XB$, where $X  = (x_{ij})$ is a $m \times n$ matrix of 
unknown variables.  Then this equation can be regarded as a 
system of linear equations with respect to the $mn$ variables 
$\{x_{ij}\}$, 
whose coefficient matrix is given by the matrix $C$.  
Moreover, if $B$ is a triangular matrix, then $C$ is block wise 
triangular, and so $\det C \ne 0$ if and only if 
$\det(A - b_{\a,\a}I_m) \ne 0$ for $1 \le \a \le n $. 
Hence we have the following lemma.
\begin{lem} 
Under the above notation, the following are equivalent.
\begin{enumerate}
\item
The matrix equation $AX = XB$ has a unique solution $X = 0$.
\item
$\det C \ne 0$.
\item
The matrices $A$ and $B$ have no common eigenvalues.
\end{enumerate}
\end{lem}
\par
We now show the following.
\begin{thm} 
Suppose that Conjecture A holds.  Then the Macdonald 
functions $P^{\pm}_{\BLa} \in \Xi_F$ are characterized by the following 
two properties.
\begin{align*}
\tag{3.6.1}
P_{\BLa}^{\pm} &= 
  m_{\Ba} + \sum_{\Bb < \Ba}w^{\pm}_{\Ba,\Bb}m_{\Bb},  \\
\tag{3.6.2}
E^{\pm}P^{\pm}_{\BLa} &= \sum_{\Bb \sim \Ba}
                   h_{\Ba,\Bb}^{\pm}P^{\pm}_{\BLa(\Bb)}.
\end{align*}
In particular, $P^{\pm}_{\BLa}$ are determined independently from 
the choice of the total order $\prec$.
\end{thm}
\begin{proof}
Let $\{ \CF_i \mid i \in I\}$ be the set of families in $Z_n^{0,0}$.
We give a total order $\prec$ on the index set $I$,  
according to the total order $\prec$ on $Z_n^{0,0}$,  and write
$H_{\pm}$ as $H_{\pm} = (H_{ij})_{i,j \in I}$, where 
$H_{ij} = H^{\pm}_{\CF_i,\CF_j}$. (Since the following discussion 
is independent of the sign $\{\pm\}$, we omit them). 
Let $X = (w^{\pm}_{\Ba,\Bb})$ be the transition matrix between basis 
$\{ m_{\Ba}\}$ and $\{P^{\pm}_{\BLa} \}$ of $\Xi_{F}$,  
and write it as $X = (X_{ij})$.
By Proposition 1.12, $X$ is block wise lower triangular, with 
identity diagonal blocks. Moreover, by Proposition 2.7, we see that
$XHX\iv = G$, where $G = (G_{ij})$ is a block
diagonal matrix with diagonal blocks $G_{ii} = H_{ii}$. 
In order to prove the theorem, we have only to show the following.
\par\medskip\noindent
(3.6.3) \ Let $X = (X_{ij})$ be a block wise lower triangular
matrix, with identity diagonal blocks, such that $XHX\iv = G$.
Then $X$ is determined uniquely, and $X_{ij} = 0$ unless 
$i = j$ or $\CF_j < \CF_i$.
\par\medskip
We show (3.6.3).  The equation $XHX\iv = G$ can be written as
\begin{equation*}
\tag{3.6.4}
H_{ii}X_{ij} - X_{ij}H_{jj} = \sum_{j \prec k \prec i}X_{ik}H_{kj}
\end{equation*}
for any pair $j \prec i$.
By backwards induction on $j$, we may assume that $X_{ik}$
are already determined for $j \prec k \prec i$. Then (3.6.4) determines
$X_{ij}$ uniquely, by Conjecture A and Lemma 3.5.  
Now suppose that $\CF_j \not < \CF_i$.  Again by induction, we may
assume that $X_{ik} = 0$ unless $\CF_k < \CF_i$. Since 
$H_{kj} = 0$ unless $\CF_j < \CF_k$ by (3.3.2), we must 
have $H_{ii}X_{ij} - X_{ij}H_{jj} = 0$.  This implies that
$X_{ij} = 0$ by Lemma 3.5, and we obtain (3.6.3).  
Thus the theorem is proved.
\end{proof}
\par\medskip\noindent
{\bf Remark 3.7.}
\addtocounter{thm}{1}
The Hall-Littlewood function $P^{\pm}_{\BLa}(x;t)$ given in [S1]
coincides with $P^{\pm}_{\BLa}(x;0,t)$.  Hence it satisfies 
similar formulas as (3.6.1), (3.6.2).  In particular, $P^{\pm}_{\BLa}(x;t)$, 
and so the Kostka functions $K^{\pm}_{\Ba,\Bb}(x;t)$ do not depend on 
the choice of the total order.  This answers the questions posed in 
[S1, Remark 4.5, (ii)] and in 
[GM, Remark 2.4], modulo the truth of the conjecture. 
\para{3.8.}
We give here some examples of the matrices $B^{\pm}_{\CF, \CF}$ and 
$H^{\pm}_{\CF,\CF}$ for some small rank cases with $e = 2$.
Assume that $W$ is the Weyl group of type $C_2$.
Then the symbols and families are given as follows.
\begin{align*}
\CF_1 &= \bigl\{ \binom {3 \ 0} {0}\bigr\}, \quad
\CF_2 = \bigl\{ \binom{2 \ 1\ 0}{2 \ 1}\bigr\}, \\
\CF_3 &= \bigl\{ \binom{2\ 0}{1}, \binom{2 \ 1}{0}, \binom{1 \ 0}{2}\bigr\},
\end{align*}
which correspond, in this order,  to the double partitions of 2,
\begin{equation*}
(2;-),\quad (-;11),\quad (1;1), \quad (11;-), \quad (-;2)
\end{equation*}
respectively.
In the tables below, we express the matrices 
$B^{\pm}_{\CF,\CF}, H^{\pm}_{\CF,\CF}$ as 
$B_{\CF}, H_{\CF}$, 
since they are independent of 
the sign $\pm$.  We have
\begin{align*}
B_{\CF_1} = (1), \quad B_{\CF_2} = (t^3q +tq), \quad 
B_{\CF_3} =
     \begin{pmatrix}
          0    & q    & tq^2 \\
          q    & 0    & -tq  \\
          tq^2 & -q^2 & 0
      \end{pmatrix}.
\end{align*}
\par
Up to $C_5$, only 3-element or 1-element families occur.
The Weyl group of type $C_6$ contains a unique 10-element family,  
which is given as follows,
\begin{equation*}
\begin{split}
\CF_4 = \bigl\{ &\binom{4\ 2\ 0}{3\ 1}, \binom{4\ 2\ 1}{3\ 0},
               \binom{4\ 3\ 0}{2\ 1},  \binom{4\ 3\ 1}{2\ 0},
               \binom{4\ 3\ 2}{1\ 0},  \\ 
               &\binom{4\ 1\ 0}{3\ 2},  \binom{3\ 2\ 0}{4\ 1},
               \binom{3\ 2\ 1}{4\ 0},  \binom{3\ 1\ 0}{4\ 2},
               \binom{2\ 1\ 0}{4\ 3}\bigr\}.
\end{split}
\end{equation*}
The corresponding double partitions of 6 are, in this order, as 
follows.
\begin{align*}
&(21;21),\quad (211;2),\quad (22;11),\quad (221;1),\quad (222;-),\\
&(2;22),\quad (11;31), \quad (111;3),\quad (1;32),\quad (-;33).  
\end{align*}
The matrix $B_{\CF}$ is given by 
\footnotesize{
\begin{equation*}
B_{\CF} =  \begin{pmatrix}
0 & q  & t^2q^3 & 0      & tq^2  & tq^2 & t^3q^4 &  0     & 0  & t^2q^3 \\
q & 0  &    0   & t^2q^3 & -tq^3 & -tq  & 0      & t^3q^4 & 0  & -t^2q^2 \\
t^2q^3 & 0 & 0 & q & -tq & -tq^3 & -t^3q^3 & 0 & t^2q^3 & 0 \\
0   & t^2q^3 & q & 0 & tq^2 & tq^2 & 0 & -t^3q^3 & -t^2q^2 & 0 \\
tq^2 & -t^2q^2 & -q^2 & tq^2 & 0 & 0 & -t^2q^2 & t^3q^2 & 0 & 0 \\
tq^2 & -q^2 & -t^2q^2 & tq^2 & 0 & 0 & 0 & 0 & t^3q^4 & -t^2q^4 \\
t^3q^4 & 0 & -t^2q^4 & 0 & -tq^3 & 0 & 0 & q & tq^2 & -t^2q^2 \\
0 & t^3q^4 & 0 & -t^2q^4 & tq^4 & 0 & q & 0 & -tq & t^2q \\
0 & 0 & t^2q^3 & -tq^3 & 0 & t^3q^4 & tq^2 & -q^2 & 0 & t^2q^3 \\
t^2q^3 & -tq^3 & 0 & 0 & 0 & -t^3q^3 & -tq^3 & q^3 & t^2q^3 & 0
          \end{pmatrix}.   
\end{equation*}
\normalsize
Throughout the above examples, $H_{\CF}$ is given by 
\begin{equation*}
  H_{\CF} = t^{-(2m+1)}B_{\CF} - 
       t^{-(2m+1)}\frac {(1 - t^m)(1 - t^{m+1})}{(1 - t)^2} 
\end{equation*}
for symbols of the shape $\Bm = (m+1, m)$, (so, $m = 1$ for
$\CF_1,\CF_3$, and $m=2$ for $\CF_2, \CF_4$, respectively).
Note that $B_{\CF}$ is not necessarily symmetric. However, if we put
$q = t$, it turns out to be symmetric since $D^1$ is a self adjoint
operator, and the representation matrix with respect to the orthonormal
basis of Schur functions coincides with the diagonal blocks
$(B_{\CF})$.
\para{3.9.}
In the case of symmetric groups, Macdonald operators are
commuting with each other since they are simultaneously diagonalizable.
In our case, Proposition 2.7 shows that $D^r_{\pm}$ are simultaneously
diagonalizable in the sense of block matrices.  
So they are commuting with each other if and only if the matrices 
$B^{r,\pm}_{\CF,\CF}$ are commuting with for $r = 1, \dots, M_1$ 
(for a fixed $\CF$).
As the following examples show, one might expect that 
$B^{r,\pm}_{\CF,\CF}$ are commuting with each other in general.
\par
First consider a simple example.  
Let $\CF_3 = \{ \BLa_1, \BLa_2, \BLa_3\}$ be 
the 3-element family of $C_2$ as in 3.8.
We consider its $m$ times shifts 
$\CF_3^{(m)} = \{ \BLa'_1, \BLa'_2, \BLa'_3\}$,
where 
\begin{equation*}
\BLa'_i = \binom { x+m,\ y+m,\ m-1, \ \dots, \ 0}
                 {z+m,\ m-1,\ \dots, \ 0} 
\end{equation*}
for $\BLa_i = \binom {x\ y}{z}$. 
Then the matrices $B^r_{\CF_3}$ are given as
\begin{equation*}
B^r_{\CF_3} = aB^1_{\CF_3} + b
\end{equation*}
with
\begin{equation*}
a = t^{2m}\sum_{J'}t^{\lp\Bdel, J'\rp},  \quad 
b = \sum_{J}t^{\lp \Bdel,J\rp}, 
\end{equation*}
where $J \in \CI_r$ runs over 
the elements of the form  
$J = \{ \Bi_{j_1}, \dots, \Bi_{j_r}\}$
for $\Bi_{j} = \binom {j+1}{j}$ with $2 \le j \le m+1$, 
and $J' \in \CI_{r-1}$ runs
over the elements having similar properties.  
This implies that $B^r_{\CF_3}$ are all commuting with for 
$1 \le r \le m$.
\par
In the following, we discuss some related results, i.e., 
we show that when $e = 2$, and $q = t$, then the operators 
$D^r$ are commuting with each other.
First we prepare a lemma.  (Since we deal with the case where
$e = 2$, we omit the sign ${\pm}$ in the discussion below.)
\begin{lem} 
Assume that $e=2$.  Then, for each $r, r'$, there exists a bijective map 
$\vf: \CI_r \times \CI_{r'} \to \CI_{r}\times \CI_{r'}$ satisfying
the following properties: let $J \in \CI_r, J'\in \CI_{r'}$ and
put $\vf(J,J') = (K, K')$.  For each $\Ba \in Z$, we have 
\begin{enumerate}
\item 
$\Ba_{JJ'} = \Ba_{K'K}$,
\item
$\lp \Ba, J\rp + \lp \Ba_J, J'\rp = 
          \lp \Ba, K'\rp + \lp \Ba_{K'}, K\rp$. 
\end{enumerate}
\end{lem}
\begin{proof}
Take $J \in \CI_{r}, J' \in \CI_{r'}$.
Let $J_0$ be the subset of $J$ consisting of 
$\Bi = \binom a b$
such that there exists 
$\Bi'= \binom x y \in J'$ with $x = a$ or 
$y = b$, and let $J_0'$ be the subset of $J'$ having similar
properties.  We define an equivalence relation on $J_0$ 
by connecting $\Bi = \binom a b, 
     \Bi' = \binom c d \in J_0$
when there exists 
$\binom x y \in J'$ such that 
$(x,y) = (a,d)$ or $(x,y) = (c,b)$.
We denote by
$\{ J_C \mid C \in \MC \}$ the set of equivalence classes in $J_0$.
Then the class $J_C$ has the following form.
\begin{align*}
J_C &= \bigl\{ \binom {a_1} {c}, \binom {a_2} {b_1}, \dots, 
           \binom {a_{k-1}}{b_{k-2}},  
               \binom {a_k} {b_{k-1}}, \binom {d}{b_k} \bigr\},
\end{align*}
where $\binom {a_i}{b_i} \in J'$ for $1 \le i \le k$.
We put 
\begin{align*}
J'_C = \bigl\{ \binom {a_i}{b_i} \mid p \le i \le q \bigr\},
\end{align*}
where $p = 0$ (resp. $q = k+1$) if there exists 
$\binom {a_0} {b_0} \in J'$ 
(resp. $\binom {a_{k+1}} {b_{k+1}} \in J'$) 
such that $b_0 = c$ (resp. $a_{k+1} = d$), and $p = 1, q= k$
otherwise.
\par
For each $C \in \MC$, we define the sets $K_C, K_C'$ as follows:
\begin{align*}
K_C = \begin{cases}
          J_C &\quad\text{ if } p=1, q = k, \\     
          J'_C
              &\quad\text{ if } p=1, q=k+1, \\
          J'_C
              &\quad\text{ if } p=0, q=k, \\
          \displaystyle\bigl\{\binom{a_0}{b_1}, \dots, 
             \binom{a_{k}}{b_{k+1}}\bigr\}
              &\quad\text{ if } p=0, q=k+1.
\end{cases}  \\
\end{align*}
\begin{align*}  
K_C' &= \begin{cases}
    \displaystyle\bigl\{ \binom{a_2}{c}, 
                 \binom{a_3}{b_1}, 
                 \dots, \binom{a_{k}}{b_{k-2}},
                         \binom{d}{b_{k-1}}\bigr\}
          &\quad\text{ if } p = 1, q = k, \\
    J_C   &\quad\text{ if } p = 1, q = k+1, \\
    J_C   &\quad\text{ if } p = 0, q = k, \\       
    J'_C
          &\quad\text{ if } p=0, q = k+1.
 \end{cases}
\end{align*}
Since $J, J'$ are subsets of the index set of elements in $Z$, 
$J_C$, etc. induce permutations on the entries of elements in $Z$.
We denote by $x_C, x'_C$ (resp. $y_C, y'_C$) the permutations in 
$\FS_{M}$  
corresponding to $J_C, J'_C$ (resp. $K_C, K'_C$), respectively.
Then it is easy to check that
\begin{equation*}
\tag{3.10.1}
x'_C\circ x_C  = y_C\circ y'_C.
\end{equation*}
We put 
\begin{equation*}
K = (J - J_0) \cup \bigcup_{C \in \MC} K_C, \quad
K' = (J' - J'_0)\cup \bigcup_{C \in \MC} K'_C.
\end{equation*}
Components of $K$ (resp. $K'$) are mutually disjoint, and we see that
$K \in \CI_{r}, K' \in \CI_{r'}$.
We now define the map 
$\vf: \CI_r \times \CI_{r'} \to \CI_r \times \CI_{r'}$
by $\vf(J,J') = (K,K')$.
Then one can check that $\vf^2 = \id$, and so $\vf$ is a 
bijection.  The assertion (i) follows from (3.10.1).
To show (ii), it is enough to verify the formula in the case
where $J = J_C, J' = J'_C, K = K_C, K' = K'_C$. 
The assertion is clear when $K_C = J'_C, K'_C = J_C$.
Assume that $p = 1, q = k$ or $p = 0, q = k+1$.  
Then one can check by a direct computation that  
\begin{equation*}
\lp \Ba, J_C\rp = \lp\Ba_{K'_C}, K_C\rp, \quad
\lp\Ba_{J_C},J'_C\rp = \lp\Ba, K'_C\rp. 
\end{equation*}
Hence the formula holds in these cases also, and the lemma
follows.
\end{proof}
\begin{prop} 
Assume that $e = 2$, and $q = t$.  Then the operators 
$D^r(t,t)$ are commuting with
each other for $r = 1, \dots, M_1$.
\end{prop}
\begin{proof}
As remarked in 3.9, it is enough to show that the matrices
$B^{r}_{\CF} = B^{r,\pm}_{\CF,\CF}$ are commuting with each 
other for $r = 1, \dots, M_1$. 
By Lemma 2.4 (i), the part corresponding to the diagonal block
in the expression of $D^r(t,t)s_{\Ba}$ is given as 
\begin{equation*}
\sum_{J \in \CI_r}t^{\lp\BLa,J\rp}s_{\BLa_J-\Bdel},
\end{equation*}
where $\BLa = \Ba +\Bdel$.
Let $\CX$ be the part of $D^{r'}(t,t)D^r(t,t)s_{\Ba}$ corresponding 
to the diagonal blocks.  Then we have
\begin{equation*}
\CX = \sum_{J \in \CI_{r}}\ve_{J}t^{\lp\BLa, J\rp}
          \sum_{J'' \in \CI_{r'}}
                   t^{\lp\, [\BLa_J],J''\rp}
                         s_{[\BLa_J]_{J''}-\Bdel}.
\end{equation*}
For each $J$, there exists $w \in \FS_{\Bm}$ such that
$[\BLa_J] = w(\BLa_J)$ and $\ve_{J}$ is given by 
$\ve_{J} = (-1)^{l(w)}$.  
Then if we put $J' =w\iv(J'')$, we have 
$\lp\, [\BLa_J], J''\rp = \lp\BLa_J, J'\rp$,
and $[\BLa_J]_{J''} = w(\BLa_{JJ'})$.
Hence one can write
\begin{equation*} 
\CX = \sum_{J \in \CI_r}\sum_{J' \in \CI_{r'}}\ve_{J,J'}
           t^{\lp\BLa,J\rp +\lp\BLa_J, J'\rp}s_{[\BLa_{JJ'}] -\Bdel},
\end{equation*}
where $\ve_{J,J'} = (-1)^{l(w')}$ for $w'\in \FS_{\Bm}$
such that $[\BLa_{JJ'}] = w'(\BLa_{JJ'})$. 
Now by applying Lemma 3.10 for $\Ba = \BLa$, we see that
$\CX$ coincides with the diagonal part for 
$D^r(t,t)D^{r'}(t,t)s_{\Ba}$.
The proposition is proved.
\end{proof}
\par\medskip

\end{document}